%
%
%
%
%
%
\input amstex
\documentstyle{amsppt}
\NoBlackBoxes

\topmatter
\title  $3x+1$ dynamics on
rationals with fixed denominator\endtitle
\author Barry Brent\endauthor

\address School of Mathematics,
University of Minnesota,
Vincent Hall,
206 Church Street S.E.,
Minneapolis, MN 55455
 \endaddress

\email brent\@math.umn.edu, barryb0\@earthlink.net 
(preferred)\endemail


\abstract We propose the existence of an 
infinite class of exact analogues of the 3x+1 conjecture 
for rational numbers 
with fixed denominators.  For some other denominators, there
 are several
attracting cycles, which exhibit scaling 
and covariance phenomena.
We analyze these phenomena 
in terms of  results of Bohm, Sontacchi and  Lagarias.

\endabstract

\thanks The author is grateful to Christian Ballot for
 helpful conversations, and to Jeff Lagarias and Keith 
Matthews for useful information provided by email.
\endthanks


\endtopmatter

\document

\head 1. Introduction\endhead   We consider dynamical 
systems on sets of
rational numbers with fixed denominator under the $3x+1$ function.  
We present numerical
evidence that there is a (probably infinite) class of these 
dynamical systems satisfying 
an exact analogue of the $3x+1$ conjecture.
We also
 discuss regularities displayed by attactors among dynamical systems of this kind that do not
happen to satisfy a $3x+1$ conjecture.

The  $3x+1$ conjecture may be stated as follows.  
Let $T: \bold{Z} \rightarrow \bold{Z}$ satisfy
$$T(x) = \cases x/2,&\text{ for $x$ even}\\
(3x+1)/2,&\text{for $x$ odd}.\endcases
$$
Let $T^n(k)$ denote the $n$th iterate of $T$ at the integer 
$k.$
The pair $(1,2)$ is an attracting cycle for $T$. 
The conjecture states that the basin
of attraction of this cycle is the entire set of positive integers.  
In other words, if $k \in \bold{Z}^+,$ then  $T^n(k) = 1$ for some  
$n$. The literature on this problem was surveyed in [Lagarias 1985] and 
[Lagarias 1998].  Generalizations are discussed in those papers 
and also in 
[Matthews 2001]. 

It is natural to  consider the effect of 
$T$ on non-positive integers.  
There have been only a few other cycles found in this way, 
namely, those containing $0, -1, -5,$
$ \text{and } -17.$  It is also natural to consider 
extensions of $T$ to  rational numbers.
The set of attracting cycles under 
such extensions 
 is the subject of 
[Lagarias 1990]. It is also discussed in [Matthews 2001]. 
Some observations 
in [Matthews 2001] are 
very much in the spirit of the present paper.
Results in [Lagarias 1990] are
used below to explain some of our own observations.
We use the notation of
[Lagarias 1990] when it is convenient.

In the sequel, fractions are understood 
to be in lowest terms unless
otherwise noted. Following Lagarias, we 
define $\bold{Q}[(2)]$ as 
the ring consisting of fractions $j/k$ with $k$ odd,
and consider $j/k \in \bold{Q}[(2)]$ to be even or 
odd according to the parity  of $j$. Then $T$
has an extension taking 
$ \bold{Q}[(2)] \text{ to } \bold{Q}[(2)]$, 
which we also denote as $T$.  

Let 
 $D_k \subset \bold{Q}[(2)]$ be the subset
 consisting of the fractions $j/k, j \in \bold{Z}^+.$ 
We say that $j/k$ has ``denominator $k$''.
The function $T$ preserves $D_k$ when 
 $k \equiv 1 \text{ or } 5 \mod 6$; 
 we restrict our attention to these cases.
The  $3x+1$ conjecture tells us that $T$ has
exactly one attracting cycle on $D_1$, the pair 
$(1,2),$ and that the basin of attraction 
of $(1,2)$ contains $D_1.$ 
(This way of putting it depends 
upon the convention, 
which we follow
throughout, that we count 
cycles that are \it irreducible \rm
in the sense of [Lagarias 1990]-- 
cycles for which the cyclic permutations are distinct.)
The behavior of $T$ on the sets $D_k$ for 
certain $k>1$ provides exact analogues
 of the original
 conjecture.
For  example:

\proclaim{Conjecture 1.1} 

\flushpar i.  The cycle of length $4$ containing $5/7$ is the 
only attractor for $D_7$, and its basin of attraction
contains $D_7$.

\flushpar ii.  The cycle of length $11$ containing $5/19$
is the only attractor for $D_{19}$, and its basin 
of attraction
contains $D_{19}$.

\flushpar iii.  The cycle of length $23$ containing $13/31$
is the only attractor for $D_{31}$, and its
 basin of attraction
contains $D_{31}$.

\endproclaim  

\flushpar The data also support,
in our opinion,

\proclaim{Conjecture 1.2} 

\flushpar
The set of positive integers $k$ such that $D_k$
lies in the basin of exactly one attractor is
infinite.

\endproclaim  
\flushpar We present numerical evidence for the two conjectures in the final section.
We will not call the following a conjecture:

\proclaim{Question  1.3} 

\flushpar
Are there divergent orbits for any of the $D_k$?
\endproclaim  
\flushpar No divergent orbits have turned up in our experiments.  On the other hand, 
 no argument  for a negative answer is known to the present writer.

The behavior of $T$ on $D_k$ for certain other 
values of $k$  
is more complex, and 
exhibits a structure.  For these $D_k$ there are 
several families of 
attractors, the members of which have equal length, 
and the lengths associated to certain of these 
families are multiples of 
the lengths associated to
 others. 
The number of odd elements exhibits 
the same scaling phenomenon as the length. Our analysis 
will explain this scaling, although some questions remain. 
For example, within the range of our observations,
usually, but not always, 
the number of odd elements of an attracting 
cycle and the length
of the cycle vary together.
This covariance 
remains enigmatic.

\head 2. Invariants of attracting cycles\endhead 
With a rational number $x \in \bold{Q}[(2)]$, we
 associate the sequence 
$$\bold{b}(x) = \{\bold{b}_j(x)\}_{j=0,1,...} 
\in \{0, 1\}^{\bold{Z}^{\geq 0}}$$ 
by setting $\bold{b}_j(x) = 0 \text{ or } 1,$ 
according to whether $T^j(x)$
is even or odd. This is the \it parity sequence 
 associated to \rm $x$.

Similarly, with any vector $\bold{c}= 
(x_0, x_1, ..., x_{n-1})$ representing a cycle 
 on 
$\bold{Q}[(2)]$, we associate a vector $\bold{v} = 
\bold{v}(\bold{c}) =
(v_0, v_1, ..., v_{n-1}) \in \{0, 1\}^n$ 
(which we term, following Lagarias,
a ``$0 - 1 $ vector'')
by setting $v_i = 0$ if $x_i$ is even, and  $v_i = 1$ if
 $x_i$ is odd. We term this the \it
parity vector associated to \rm  $\bold{c}.$
This association is one-to-one.
The following 
uniqueness theorem
(according to Lagarias, ``essentially due to B\"ohm and 
Sontacchi'' [B\"ohm and Sontacchi 1978])
appears in [Lagarias 1990] (p. 36).  We will
refer to it as the 

\proclaim{ BSL Theorem} Given any $0 - 1$ vector 
$\bold{v} = (v_0, v_1, ..., v_{n-1})$
there is a unique $x$ in  $\bold{Q}[(2)]$ which is
 periodic of period $n$ under iteration
by the $3x+1$ function $T$, and whose parity sequence 
starts with $\bold{v}$.  It is
given by
$$ x = x(\bold{v}) = (2^n - 3^{v_0 + ... +v_{n-1}})^{-1} 
\sum_{j=0}^{n-1} v_j 3^{v_{j+1}+...+v_{n-1}}2^j.
\tag{2.1}
$$

\endproclaim

Evidently, $n$ is the length of  $\bold{v}$. 
We introduce the notation $\lambda(\bold{v}) := n.$
We introduce two more invariants: 
$\omega(\bold{v}) :=\sum_{j=0}^{n-1} v_j,$ and 
 $\rho(\bold{v}) := \sum_{j=0}^{n-1} v_j 3^{v_{j+1}+...+v_{n-1}}2^j.$ 
Thus (suppressing the dependence on
$\bold{v}$) we have
$$ x =     \frac {\rho}{2^ {\lambda}- 3^{\omega}}.
\tag{2.2}
$$

\head 3. The number of attractors for a given denominator\endhead

Let $\alpha(k)$ be the  number of attractors for 
$T$ on $D_k$, 
and let $\alpha_n(k)$ be the number of such attractors of length $n$.
The BSL Theorem allows an interpretation of
the numbers
 $\alpha(k)$ and  $\alpha_n(k)$. 
Let $V(k)$ be the set  
of $0 - 1$ vectors $\bold{v}$ such that
$k = 
\text{ the denominator of } x(\bold{v})$,
\it i.e. \rm
 $k = (2^ {\lambda}- 3^{\omega})/(\rho, 2^ {\lambda}- 3^{\omega}).$
Let $\nu(k)$ be the size of $V(k)$. Let $V_n(k)$ be the subset 
of $V(k)$  with
$\lambda=n,$  and let $\nu_n(k)$ be the size of $V_n(k).$
Then

\proclaim{ Proposition 3.1} 

1. 
$\alpha(k)$ is finite if 
and only if  $\nu(k)$
is finite.

2.
$\alpha(k) \neq 0$  if and only if $\nu(k) \neq 0$.

3.  $\alpha_n(k)=\nu_n(k)/n.$

4. $\alpha(k) = \sum_{n=1}^{\infty} \nu_n(k)/n .$ 
\endproclaim

\flushpar \it \,\,Proof \rm Only the third claim needs an argument.
\,\, If  $\bold{v} \in V_n(k),$ the cyclic permutations of  $\bold{v}$
are distinct.
Otherwise, distinct cyclic permutations of the attractor 
containing $x(\bold{v})$
correspond to non-distinct cyclic permutations of $\bold{v}$, and
this is ruled out by the BSL Theorem. Furthermore,
the cyclic permutations of $\bold{v}$ 
all lie in $V_n(k)$, because they correspond 
to cyclic permutations of the attractor starting with $x(\bold{v})$,
and,
if $k \equiv 1 \text{ or } 5 \mod 6$, 
all numbers in the cycle 
corresponding to
$\bold{v}$ have denominator $k$. \boxed{}
\vskip .1in
We can use Proposition 3.1 and a result of Lagarias
to derive an identity for  $\nu_n$.  Let $I(n)$ count the
number of irreducible cycles of length $n$. It is easy to see that
the denominators $k$ of cycles in $\bold{Q}[(2)]$ must all satisfy
$k \equiv 1 \text{ or } 5 \mod 6$, so 
$I(n) = \sum_{k=1}^{\infty} \alpha_n(k) = 
\sum_{k=1}^{\infty} \nu_n(k)/n.$ 
Now, 
 Eq. 2.5 of [Lagarias 1990] (p. 38) states:

$$
I(n) = \frac 1n \sum_{d|n} \mu(d) 2^{n/d} .
$$
By combining these facts, we arrive at

\proclaim{ Proposition 3.2}
 
$$\sum_{k=1}^{\infty} \nu_n(k) =  \sum_{d|n} \mu(d) 2^{n/d}.$$

\endproclaim

\flushpar The quantity $\nu_n(k)$, $n$ fixed, vanishes 
for large $k$ (obvious). It would
be interesting to 
have a good bound on  $k$ to replace the upper bound in 
the left-hand sum.

\head 4. Scaling, covariance and repetition\endhead
The following table illustrates typical behavior of the parameters we studied. 
The denominator
is specified by the parameter $k$. 
We specify
attractors $\bold{c}$ by the numerators of their smallest members.
For example, the first row conveys the following: there is a 
cycle for $D_5$ with  smallest
member $ 1/5$, length three and one odd element. 
Since the parameters $\lambda(\bold{v})$ and $\omega(\bold{v})$ 
are invariant under cyclic permutation of $\bold{v},$ they are, in effect,
invariants of $\bold{c}(\bold{v})$, namely, the length of 
 $\bold{c}$ and the number of odd elements of  $\bold{c}$,
respectively.
 We treat them as such in this table and later. 

\flushpar 
\settabs 6 \columns
\+\cr\+&&$k$ & $\bold{c}$ & $\lambda(\bold{c})$ & $\omega(\bold{c})$\cr
\+\cr\+&&$5$ & 1 & 3 & 1\cr \vskip -.15in
\+\cr\+&&$5$ & 19 & 5 & 3\cr \vskip -.15in
\+\cr\+&&$5$ & 23 & 5 & 3\cr \vskip -.15in
\+\cr\+&&$5$ & 187 & 27 & 17\cr \vskip -.15in
\+\cr\+&&$5$ & 347 & 27 & 17\cr \vskip -.15in
\+\cr\+&&$7$ & 5 & 4 & 2\cr \vskip -.15in
\+\cr\+&&$11$ & 1 & 6 & 2\cr \vskip -.15in
\+\cr\+&&$11$ & 13 & 14 & 8\cr \vskip -.15in
\+\cr\+&&$13$ & 1 & 4 & 1\cr \vskip -.15in
\+\cr\+&&$13$ & 131 & 24 & 15\cr \vskip -.15in
\+\cr\+&&$13$ & 211 & 8 & 5\cr \vskip -.15in
\+\cr\+&&$13$ & 227 & 8 & 5\cr \vskip -.15in
\+\cr\+&&$13$ & 259 & 8 & 5\cr \vskip -.15in
\+\cr\+&&$13$ & 251 & 8 & 5\cr \vskip -.15in
\+\cr\+&&$13$ & 283 & 8 & 5\cr \vskip -.15in
\+\cr\+&&$13$ & 287 & 8 & 5\cr \vskip -.15in
\+\cr\+&&$13$ & 319 & 8 & 5\cr
\+\cr\flushpar

\flushpar We will refer to a search among the fractions
$j/k, 1 \leq j \leq N$ as a search of depth $N$.
For the given $k$, the table lists all attractors $\bold{c}$ for $D_k$ 
such that 
some of the fractions surveyed in a search of
depth $500$ fall
into  $\bold{c}$.

The $k=13$ rows illustrate the phenomenon 
we term \it scaling: \rm (1)
for a given $D_k$ and a more  or less numerous set of attractors
of equal length,
 the ratios $\lambda/\omega$ are identical
 and are equal
to the corresponding ratio for a single 
additional attractor of greater length. 
(2) For a given denominator, when two attractors
$\bold{c}_1$ and $\bold{c}_2$ satisfy 
 $\lambda(\bold{c}_1)/\omega(\bold{c}_1)=\lambda(\bold{c}_2)/\omega(\bold{c}_2)$, 
it is also the case that  $\lambda(\bold{c}_1)/\lambda(\bold{c}_2)$ and 
$\omega(\bold{c}_1)/\omega(\bold{c}_2)$ are integers.
The $k=5,11$ rows exhibit dynamical systems without scaling.
Of course it is possible that sufficiently deep searches 
would turn up scaling for every $k$.
This would refute both Conjectures 1 and 2. 
\it
Covariance \rm is our term for the apparent interdependence 
for a fixed $D_k$ of
$\omega(\bold{c})$ and $\lambda(\bold{c})$: they seem to 
vary together. There are  
exceptions to this behavior, but its frequent occurrence
seems to require an explanation.  
The smallest examples in which covariance fails are given in the
following table.

\flushpar 
\settabs 6 \columns
 \+\cr\+&&$k$ & $\bold{c}$ & $\lambda(\bold{c})$ & $\omega(\bold{c})$\cr
\+\cr\+&&$511$ & 11 & 54 & 24\cr \vskip -.15in
\+\cr\+&&$511$ & 293 & 45 & 24\cr \vskip -.15in
\+\cr\+&&$757$ & 43 & 84 & 37\cr \vskip -.15in
\+\cr\+&&$757$ & 85 & 84 & 46\cr \vskip -.15in
\+\cr
\vskip .1in \flushpar
 Finally, we see in our data 
\it repetition \rm of identical $\lambda$ and $\omega$
values for distinct attractors of a single denominator.  Our 
first table supplies several examples.  In a search of depth $50$ among
all denominators $k \equiv 1 \text{ or } 5 \mod 6, 
1 \leq k \leq 1501,$ we found 121 examples of scaling and 205 examples
of repetition.  There were 83 denominators exhibiting both of these behaviors.

Some of this we can explain, but not all of it.  
Even if our conjectures are correct, we can ask 
whether
there are $D_k$ with more than one attractor
that do not exibit scaling.  We do not know
the answer.
Below, we explain 
the existence of scaling, but  not
 the pattern in which a single long attractor shares
the same  $\lambda/\omega$ ratio as a larger set of short attractors.
We cannot explain covariance.  We offer speculation instead.  
From the BSL Theorem, it is clear that
the parity vectors for different $\bold{c}$ are never  
related by cyclic permutations, yet some families
of attracting cycles do appear to be related by 
non-cyclic permutations of their parity vectors; 
superficially
this is what covariance comes down to.
We can ask: do these permutations form a group or some other structure?
A positive answer might create a situation in which
 we could predict the
existence of cycles missed by our searches 
(because of their insufficient depth) to 
fill out empty slots in such a structure.

Next, we give an account of scaling and repetition that covers all 
the examples of these phenomena in our data.
We cannot rule out the existence of sporadic examples that come about some other way. 
Let $k | J \text{ and } J = dk.$
To a representation of $J$ in the form 
$ J = 2^{\lambda} - 3^{\omega},$ 
we associate  a 
possibly empty  set $V_{\lambda, \omega,  d} $ of
$0 - 1$ vectors $\bold{v}$ with the properties that
$\lambda(\bold{v}) = \lambda, \omega(\bold{v}) = \omega,$ and
$(\rho(\bold{v}), J) = d.  $ 
For these $\bold{v},$ the BSL Theorem provides  
 an  attractor $ \bold{c}(\bold{v})$ 
with 
$\lambda(\bold{c}(\bold{v})) = \lambda$ and 
$\omega(\bold{c}(\bold{v})) = \omega$,
containing $$x(\bold{v}) = \frac {\rho(\bold{v})}{J} \in D_{J/d} = D_k.$$

If it happens that there is a pair $\bold{v}_1, \bold{v}_2$ in $ V_{\lambda, \omega,  d}$
that
are not related by a cyclic permutation, 
then  $\bold{c}(\bold{v}_1) \neq \bold{c}(\bold{v}_2)$ but
$\lambda(\bold{c}(\bold{v}_1)) = \lambda(\bold{c}(\bold{v}_2))$ and 
$\omega(\bold{c}(\bold{v}_1)) = \omega(\bold{c}(\bold{v}_2))$, 
which is repetition.

Now suppose that 
$ s =\lambda \delta, t = \omega \delta, \delta \in \bold{Z}^+ .$ 
Then  $k | 2^s - 3^t = M =fk$ (say).
If   $V_{\lambda, \omega,  d} $ is
 non-empty  then fractions $x(\bold{v})$  and attractors
$\bold{c}(\bold{v})$ exist satisfying the conditions above. 
If also some $\bold{w} \in V_{s, t,  f} $, 
then so do
fractions $x(\bold{w})$  and attractors
$\bold{c}(\bold{w})$ with 
$x(\bold{w}) \in \bold{c}(\bold{w})$,
$\lambda(\bold{c}(\bold{w})) = s,$  
$\omega(\bold{c}(\bold{w})) = t$, and
 $$x(\bold{w}) = \frac {\rho(\bold{w})}{M} \in D_{M/f} = D_k.$$
Scaling arises because 
$ s/\lambda = t/\omega.$  The integrality we mentioned above arises because 
$\delta$ is an integer. One may ask whether the integrality in 
turn reflects unknown
properties of the parity vectors.

\head 5. Some evidence for Conjectures 1.1 and 1.2 \endhead

The evidence we can offer for Conjecture 1.1 is straightforward.  We executed a 
depth $3 \times 10^6$ search 
for denominator $7$, a depth $1.3 \times10^6$  search
for denominator $19$ 
and a depth $ 10^6$ search for denominator $31$ 
without turning up a second attracting cycle or a divergent orbit. 

To test Conjecture 1.2, we searched for denominators $k, 1 \leq k \leq 2000,$
with a single attracting cycle. We conducted
searches of varying depth and examined the effect of this variation.  
In the following table,
 we display  the number $A$ of denominators in 
the given range showing just one attracting cycle
as a function of the search depth $N$.

\flushpar 
\settabs 6 \columns
 \+\cr\+&&$N$ & $A$ &$N$ & $A$ \cr
\+\cr\+&&$20$ & 213 &$800$ & 166\cr \vskip -.15in
\+\cr\+&&$50$ & 184 &$1600$ & 166\cr \vskip -.15in
\+\cr\+&&$100$ & 181 &$2400$ & 166\cr \vskip -.15in
\+\cr\+&&$200$ & 176 &$3200$ & 162\cr \vskip -.15in
\+\cr\+&&$400$ & 172\cr \vskip -.15in

\+\cr
\vskip .1in \flushpar

The data is not obviously inconsistent with a model of the 
form $A=c_1 + f(N), 0 \leq f(N) \leq (2000-c_1) \exp(-c_2 N),$ for constants
$c_1, c_2.$ Of course, this model, if it
is correct, tends to confirm Conjecture 1.2.  
The \it Mathematica \rm Statistics package \it NonlinearFit \rm
command picks out the values 
$c_1 = 171.594, c_2 = 0.189263$ for this data set and the cruder model
$A=c_1 + (2000-c_1) \exp(-c_2 N)$, but it is
insensitive to large perturbations of the data, and so we need not take the
particular choice seriously.  However, \it Mathematica \rm refrained from issuing
a warning message, not a bad sign, at least, and  the only measure 
we have of the goodness of the model's
fit.
\refstyle{C}

\enddocument